\documentclass{article}
\usepackage{amssymb,amsmath,amsthm}
\newtheorem{theorem}{Theorem}

\newtheorem{lemma}{Lemma}

\newtheorem{remark}{Remark}
\date{}
\numberwithin{equation}{section} \numberwithin{theorem}{section}
\numberwithin{lemma}{section} \numberwithin{corollary}{section}
\numberwithin{remark}{section} \numberwithin{proposition}{section}
\numberwithin{definition}{section}
\begin{document}
\newcommand{\n}{\noindent}
\newcommand{\vs}{\vskip}

\title{ Lipschitz Continuity of Solutions to a Free Boundary Problem involving the $A$-Laplacian}

\author{S. Challal$^1$ and A. Lyaghfouri$^{2}$\\
\\
$^1$ Glendon College, York University\\
Department of Mathematics\\
Toronto, Ontario, Canada\\
\\
$^2$ American University of Ras Al Khaimah\\
Department of Mathematics and Natural Sciences\\
Ras Al Khaimah, UAE\\
}\maketitle

\vs 1cm
\begin{abstract}
In this paper we prove local interior and boundary Lipschitz continuity
of solutions of a free boundary problem involving the $A$-Laplacian.
We also show that the free boundary is represented
locally by graphs of a family of lower semi-continuous functions.
\end{abstract}

\vs 0.5cm

\n Key words : $A$-Laplacian, Free boundary, Lipschitz continuity

\vs 0.5cm \n AMS 2000 Mathematics Subject Classification: 35R35; 35J60

 \vs 0.5cm
\section*{Introduction}\label{S:intro}

\vs 0.5cm \n We consider the following problem

\begin{equation*}(P)
\begin{cases}
& \text{ Find }\,\,(u, \chi) \in  W^{1,A}(\Omega)\times L^\infty
(\Omega)
\text{ such that } :\\
& (i)\quad  0\leqslant u\leqslant M, \quad 0\leqslant
\chi\leqslant 1 ,
\quad u(1-\chi) = 0 \,\,\text{ a.e.  in } \Omega\\
 & (ii)\quad \Delta_A u=-div(\chi H(x)) \quad\text{in }
 (W_0^{1,A}(\Omega))',
\end{cases}
\end{equation*}

\n where $\Omega$ is an open bounded domain of
$\mathbb{R}^n$, $n\geqslant 2$, $x=(x_1,...,x_n)$, $M$ is a positive constant,
 $H(x)$ is a vector function and  $\Delta_A$ is the $A$-Laplacian given by

$$\Delta_A u=div\Big({{a(|\nabla u|)}\over{|\nabla u|}}\nabla
u\Big)$$

\n with $\displaystyle{A(t)=\int_0^t a(s)ds}$, $a$ is  a  $C^1$ function from
$ [0,\infty)$ to $ [0,\infty)$ such that $a(0)=0$ and satisfies for some positive
constants $a_0, a_1$

\begin{equation}\label{0.1}
a_0\leqslant {{ta'(t)}\over{a(t)}}\leqslant a_1    \qquad \forall t>0.
\end{equation}

\n As a consequence of (0.1), we have the following
monotonicity inequality that was established in \cite{[CL3]}

\begin{equation}\label{0.2}
 \Big( { a(|\xi|)\over |\xi|}  \xi-  { a(|\zeta|)\over |\zeta|}
 \zeta\Big). (\xi-\zeta) >0 \qquad \forall \xi, \zeta\in
 \mathbb{R}^n\setminus\{ 0\},\quad \xi\neq\zeta.
\end{equation}

\n We point out that there is a wide range of such functions. To
start with, observe that we have $a_0=a_1=p-1$ if and only if
$a(t)=t^{a_0}$. In this case we have $A(t)={{t^p}\over p}$ and
$\Delta_A$ is the $p-$Laplacian operator $\Delta_p$. A second class
of such functions is given by
$$a(t)=\left\{
\begin{array}{ll}
    c_1t^{\alpha}, & \hbox{if } 0\leq t<t_0 \\
    c_2t^{\beta}+c_3, & \hbox{if } t\geqslant t_0,\\
\end{array}
\right.
$$
\n where $t_0$, $\alpha$ and $\beta$ are positive numbers, $c_1$,
$c_2$ and $c_3$ are real numbers such that $a(t)$ is a $C^1$
function. In this case (0.1) is satisfied with
$a_0=\min(\alpha,\beta)$ and $a_1=\max(\alpha,\beta)$.

\vs 0.2cm \n A third class of these functions is given by
$a(t)=t^\alpha\ln(\beta t+\gamma)$, with $\alpha, \beta,
\gamma>0$, in which case (0.1) is satisfied with $a_0=\alpha$ and
$a_1=1+\alpha$.

\vs 0,5cm\n We assume that the vector function $H=(H_1,...,H_n)$
satisfies for some positive constant $\overline{h}$
\begin{eqnarray}\label{0.2-3}
 & |H|_\infty\leqslant \overline{h},\\
&|div(H)|_\infty\leqslant \overline{h}.
\end{eqnarray}

\vs 0.2cm \n The natural framework for the $A$-Laplacian is the Orlicz-Sobolev space
$W^{1,A}(\Omega)$. Here we recall the definition of those spaces with their respective
norms under which they are Banach, reflexive vector spaces:

\begin{eqnarray*}
&& L^A(\Omega)=\left\{\,u\in L^1(\Omega)\,:\,\int_\Omega
A(|u(x)|)dx<\infty\,\right\},
~~||u||_A=\inf\left\{k>0\,:\,\int_\Omega A\Big({{|u(x)|}\over k}\Big)dx\leqslant 1\,\right\}\\
&& W^{1,A}(\Omega)=\left\{\,u\in L^A(\Omega)\,:\,|\nabla u|\in
L^A(\Omega)\,\right\},\qquad  ||u||_{1,A}=||u||_A+||\nabla u||_A.
\end{eqnarray*}

\n $L^A(\Omega)$ and $ W^{1,A}(\Omega)$ are  Banach, reflexive spaces.

\vs 0.2cm\n Throughout this paper, we shall denote by $B_r(x)$ a ball with center
$x$ and radius $r$, and by $B_1$ the unit ball.

\vs 0.2cm\n we recall that problem $(P)$ in the case of a linear operator was considered in
\cite{[C]}, \cite{[CL3]}, \cite{[CL5]}, and \cite{[Ly1]}. The case of the
$p-$Laplacian was considered in \cite{[Cha]}, \cite{[CaL]}, \cite{[CL1]},
\cite{[CL2]}, \cite{[CL6]}, and \cite{[CL7]}.
Regarding the problem with a Newman boundary condition, we refer for example to
\cite{[ChiL1]}, \cite{[ChiL2]}, \cite{[Ly2]}, \cite{[Ly3]}, and \cite{[Ly4]}.

\vs 0.5cm\n According to a result established in \cite{[CL8]}, we
know that $u\in C_{loc}^{0,\alpha}(\Omega)$ for some $\alpha\in(0,1)$
provided that $H\in L^\infty_{loc}(\Omega)$.
In this paper, we will show that $u\in C_{loc}^{0,1}(\Omega\cup T)$, where $T$ is a
nonempty open connected subset of $\partial\Omega$ on which $u$ satisfies the
Dirichlet condition $u=0$. In the second part of the paper, we assume that $div H
\geqslant 0$ and we establish that the free boundary is represented by
a family of lower semi-continuous functions.

\n This work generalizes results from  \cite{[CL1]}, \cite{[CL5]}, and \cite{[CL6]}
for the p-Laplacian and results in \cite{[C]}, \cite{[CL2]}, and \cite{[CL4]},
for a linear operator. To avoid repetitions, we shall refer the reader to those
works for any omitted technical details.

\vs 0,5cm

\section{Interior and boundary Lipschitz continuity}\label{1}

\vs 0,5cm\n The first main result of this section is the following interior regularity.

\vs 0,3cm\n
\begin{theorem}\label{t1.1} Let $(u, \chi)$ be a solution of
$(P)$. Then $ u\in C^{0,1}_{loc}(\Omega).$

\end{theorem}

\vs 0,3cm \n We observe that since $H\in L^\infty_{loc}(\Omega)$, we
have $u\in C_{loc}^{0,\alpha}(\Omega)$ for some $\alpha\in(0,1)$ \cite{[CL8]}.
Consequently the set $\{u>0\}$ is open. Moreover we have $\Delta_A u=-div(H)$ in ${\cal
D}'(\{u>0\})$ and $div(H)$ is uniformly bounded in $\Omega$. So
we have $ u\in C^{1,\beta}_{loc}(\{u>0\})$ for some $\beta\in(0,1)$
\cite{[L]}. Therefore to prove Theorem 1.1, it is enough
to study the behavior of $u$ near the free boundary. This is the object
of the following lemma.

\vs 0,5cm

\begin{lemma}\label{l2.1} Let $x_0=(x_{01},...,x_{0n})$ and $r>0$ such that
$B_r(x_0)\subset \{u>0\}$, $\overline{B_r(x_0)}\subset\Omega$
and $\partial B_r(x_0)\cap \partial \{u>0\}\neq \emptyset$. Then there
exists a positive constant $C$ depending only on $n$,
$\overline{h}$, $a_0$, $a^{-1}(\overline{h})$, and $\delta(\Omega)$
(the diameter of $\Omega$) such that
$$\sup_{B_{r/2}(x_0)}u \, \leqslant C \, r.$$

\end{lemma}

\vs 0,5cm

\n\emph{Proof.} We start by applying Harnack's inequality (see \cite{[L]}, Corollary 1.4):
$$\sup_{B_{r/2}(x_0)}u\,\leqslant \, C\big(  \inf_{B_{r/2}(x_0)}u
    \,+\,r. a^{-1}(\overline{h}\delta(\Omega))\big),$$
where $C$ is a positive constant depending only on $n$, $a_0$ and
$a_1$.

\n Therefore, to prove the lemma, it will be enough to establish the
inequality
$$\displaystyle{\min_{\overline{ B_{r/2}(x_0)}} u \,\leqslant \,C
\,r}.$$

\n So let $\epsilon \in(0,r)$ such that
$B_{r+\epsilon}(x_0)\subset \Omega$ and let $v$ be defined in the circular ring
$D=B_{r+\epsilon}(x_0)\setminus \overline{B_{r/2}(x_0)}$ by
\[v(x) = k\big(
 e^{-\alpha\rho^{2}}-e^{-\alpha(r+\epsilon)^{2}}\big)\]
where
\[\rho= |x-x_0|,
  \quad k = \displaystyle{{m\over{e^{-\alpha
    r^2/4}-e^{-\alpha(r+\epsilon)^2}}}},
    \quad \displaystyle{m=\min_{\overline{
B_{r/2}(x_0)}} u}, \quad \alpha=\displaystyle{{\kappa\over{r^2}}},
\quad \displaystyle{\kappa=2\Big(1+{n\over{a_0}}\Big)}.\]

\vs 0,2cm\n We claim that

\begin{equation}\label{e1.1}
  \Delta_A  v\geqslant {a(|\nabla v|) \over\rho}\qquad \hbox{ in }  D.
\end{equation}

\n Indeed, we first observe that
\begin{equation}\label{e1.2}
 \displaystyle{ \Delta_A v= {{a(|\nabla v|)}\over{|\nabla v|^3}}\Big\{|\nabla v|^2
  \Delta v+\Big({{a'(|\nabla v|)}\over{a(|\nabla v|)}}|\nabla
  v|-1\Big)\sum_{i,j}{{\partial v}\over{\partial x_i}}{{\partial v}\over{\partial x_j}}
   {{\partial^2 v}\over{\partial
x_i\partial x_j}}\Big\}}.
\end{equation}

\n Moreover we have in $D$

\begin{eqnarray*}
&& \nabla v=-2\alpha k e^{-\alpha \rho^{2}}(x-x_0), \qquad |\nabla
v|=2\alpha k \rho e^{-\alpha \rho^{2}},
 \qquad \Delta v=-2\alpha k e^{-\alpha \rho^{2}}(n-2\alpha \rho^2),\\
&& {{\partial^2 v}\over{\partial x_i\partial x_j}}=-2\alpha k
e^{-\alpha \rho^{2}}\Big(\delta_{ij}-2\alpha(x_i-x_{0i})(x_j-x_{0j})\Big),\\
&& \sum_{i,j}{{\partial v}\over{\partial x_i}}{{\partial
v}\over{\partial x_j}} {{\partial^2 v}\over{\partial x_i\partial
x_j}} =-(2\alpha k )^3  \rho^2e^{-3\alpha
\rho^2}(1-2\alpha\rho^2).
\end{eqnarray*}

\n Taking into account the fact that $\displaystyle{
1-2\alpha\rho^2=1-2{\kappa\over r^2}\rho^2\leqslant
1-2{\kappa\over r^2}\Big({r\over 2}\Big)^2=1-{\kappa\over 2}<0}$,
we get by substituting the above formulas in (1.2)

\begin{eqnarray*}
&& \displaystyle{ \Delta_A v= -(2\alpha k)^3
\rho^2e^{-3\alpha\rho^2}{{a(|\nabla v|)}\over{|\nabla
 v|^3}}\Big\{n-1
 +{{a'(|\nabla v|)}\over{a(|\nabla v|)}}|\nabla
  v|(1-2\alpha\rho^2)\Big\}}\nonumber\\
 && \qquad \geqslant -(2\alpha k)^3
\rho^2e^{-3\alpha\rho^2}{{a(|\nabla v|)}\over{|\nabla
 v|^3}}\Big\{n-1
 +a_0(1-{\kappa\over 2})\Big\}\qquad\hbox{ by } (0.1)\nonumber\\
 &&\qquad  =-{{a(|\nabla v|)}\over{\rho}}\Big(n-1
 +a_0(1-{\kappa\over 2})\Big)={{a(|\nabla v|)}\over{\rho}}.
\end{eqnarray*}

\n Hence (1.1) holds.

\vs 0.2cm\n Then we deduce from (0.4) and (1.1) that

\begin{eqnarray}\label{e1.3}
&&\Delta_A v+div(H)\geqslant {{a(|\nabla
v|)}\over{\rho}}-\overline{h}\nonumber\\
&& \displaystyle{\geqslant {1\over {(r+\epsilon)}}a\Big(
2{\kappa\over r^2}.{{me^{-{\kappa\over r^2}(r+\epsilon)^2}}\over
 {e^{-{\kappa\over 4}}-e^{-{\kappa\over r^2}(r+\epsilon)^2}}}.{r\over
 2}\Big)-\overline{h}=\theta(r).}
\end{eqnarray}

\n * If  $\quad$  $\theta(r)\leqslant 0$, $\quad$ then
$\quad$ $\displaystyle{ a\Big( {\kappa\over r}.{{me^{-{\kappa\over
r^2}(r+\epsilon)^2}}\over
 {e^{-{\kappa\over 4}}-e^{-{\kappa\over r^2}(r+\epsilon)^2}}}\Big)\leqslant \overline{h}(r+\epsilon)}$.
 $\quad$   Letting  $\epsilon\rightarrow 0$,  we get $\quad$
$ \displaystyle{a\Big( {\kappa\over r}.{{me^{-\kappa}\over
 {e^{-{\kappa\over 4}}-e^{-\kappa}}}}\Big)\leqslant \overline{h}r\leqslant \overline{h}\delta(\Omega)}$.
Hence

$$m\leqslant
a^{-1}\big(\overline{h}\delta(\Omega)\big) {{(e^{{3\over 4}\kappa
}-1)}\over \kappa}r =C(n,a,\overline{h},\delta(\Omega))r$$
\n
 and the lemma follows.

\vs 0.2cm\n * If $\quad$  $\theta(r)>0$,  $\quad$ then we deduce from  (1.3),
since $\pm(v-u)^+ \in W^{1,A}_0(D)$, that

\begin{equation}\label{1.4}
 \displaystyle{\int_D }\Big( {{a(|\nabla v|)}\over{|\nabla v|}}\nabla
v + H(x)\Big)  .\nabla(v-u)^+ \,\leqslant\, 0.
\end{equation}

\vs 0.2cm\n From $(P)ii)$, we also have

\begin{equation}\label{1.5}
\displaystyle{\int_D }\Big( {{a(|\nabla u|)}\over{|\nabla u|}}\nabla
u + \chi H(x)\Big)  .\nabla(v-u)^+ \,=\, 0.
\end{equation}

\n Subtracting (1.5) from (1.4), we get

\begin{equation*}
 \int_{D} \Big( {{a(|\nabla v|)}\over{|\nabla v|}}\nabla
v-{{a(|\nabla u|)}\over{|\nabla u|}}\nabla u \Big).\nabla(v-u)^+
dx \leqslant
 \int_{D}(\chi-1)H(x).\nabla(v-u)^+ dx
\end{equation*}

\n which can be written  by using (0.3) and the
fact that $\chi=1$ a.e. in $\{u>0\}$

\begin{eqnarray*}
&& \int_{D\cap \{u>0\}} \Big( {{a(|\nabla v|)}\over{|\nabla v|}}\nabla
v-{{a(|\nabla u|)}\over{|\nabla u|}}\nabla u \Big).\nabla(v-u)^+
dx\\
&&\leqslant  \int_{D\cap [u=0]}(\chi-1)H(x).\nabla v dx
-\int_{D\cap [u=0]} {{a(|\nabla v|)}\over{|\nabla
v|}}|\nabla v|^2 dx\\
&& \leqslant \int_{D\cap [u=0]} |\nabla
  v| \big(\overline{h}- a(|\nabla v|)\big)dx.
\end{eqnarray*}

\vs 0.2 cm\n If $ \displaystyle{ \int_{D\cap [u=0]}
|\nabla v| \big(\overline{h}- |a(|\nabla v|)|\big)dx\leqslant 0}$, then
we would get by taking into account (0.2) that $\nabla (v-u)^+=0$ in $B_r(x_0)$. Since $(v-u)^+=0$
on $\partial B_r(x_0)$, this leads to $v\leqslant u$ in $\overline{D\cap B_r(x_0)}$.
But given that $v>0$ on $\partial B_r(x_0)$ and $\partial B_r(x_0)\cap\partial
\{u>0\}\neq\emptyset$, we would get a contradiction. Hence

\begin{equation}\label{e1.6}
\displaystyle{ \int_{D\cap [u=0]} |\nabla
  v| \big(\overline{h}- a(|\nabla v|)\big)dx>0}.
\end{equation}

\n Since $|\nabla v|= 2k\alpha\rho e^{-\alpha\rho^{2}}$ and
$\kappa\geqslant 2$, we have
$${d\over{d\rho}}
|\nabla v|= 2k\alpha e^{-\alpha\rho^{2}}(1-2\kappa{\rho^2\over
r^2}) \leqslant 2k\alpha
e^{-\alpha\rho^{2}}(1-\kappa/2)\leqslant 0.$$

\n Therefore $|\nabla v|$ is non-increasing with respect to $\rho$.
It follows  from (1.6) that
 $\displaystyle{{a(|\nabla v|)}_{|_{\partial B_{r+\epsilon}(x_0)}}}$ $\displaystyle{=a( 2k\alpha (r+\epsilon)
 e^{-\alpha
 (r+\epsilon)^2})< \overline{h}}$
   i.e.
 $\displaystyle{ {{2m\kappa(r+\epsilon)e^{-\alpha(r+\epsilon)^2}}\over{r^2\big(e^{-\alpha
r^2/4}-e^{-\alpha(r+\epsilon)^2}\big)}}< a^{-1}(\overline{h})}$.
Letting $\epsilon \rightarrow 0$, we obtain
$  \displaystyle{{{2m\kappa
e^{-\kappa}}\over{r\big(e^{-\kappa/4}-e^{-\kappa}\big)}}\leqslant
a^{-1}(\overline{h})}$   which leads to

$$m\leqslant {{a^{-1}(\overline{h})} \over {2\kappa}} (e^{3\kappa/4}-1) r= C(\overline{h},\kappa,a)\, r.$$ \qed

\vs 0,5cm \n \emph{Proof\,\, of \,\,Theorem 1.1.} The
proof is based on Lemma 1.1 and arguments similar to those in
the p-Laplacian case \cite{[CL6]}.  In particular, we use the
scaling function
\[v(y)= {{u(x_0+Ry)}\over R}\quad\text{for }\quad y\in B_1\]

\n which satisfy the equation

$$\Delta_A v =-R(div H)( x_0 + R y)\quad\text{in}\quad B_1\quad\text{ if
}\quad B_R(x_0) \subset\{u>0\}$$

\n and to which we apply the estimate from \cite{[L]}, Theorem 1.7. So we get for
some positive constant $C(n,a,M,R)$
\[\sup_{B_{1/2}}|\nabla v|\,\leqslant \,C(n,a,M,R).\]
\qed

\vs 1cm Now we assume that $u=0$ on a nonempty subset $T$ of $\partial\Omega$,
and we study the Lipschitz continuity of $u$ up to $T$.
We assume that the uniform exterior sphere condition is satisfied
locally on $T$ i.e. for each open and connected subset
$S_0\subset\subset T$

\begin{equation}\label{1.9}
\exists R_0 >0\quad\text{such that}\quad\forall x\in S_0
\quad\exists y_0\in \mathbb{R}^n\setminus \overline{\Omega}\quad
\overline{B_{R_0}(y_0)}\cap S_0=\{x_0\}.
\end{equation}

\n We observe that we can assume that
$R_0<d_0/3\quad\text{where}\quad d_0=d(S_0,\partial\Omega\setminus
T)>0$, which we will do throughout this section.
Then we have our second main result:

\vs 0,5cm
\begin{theorem}\label{t1.2} For any solution $(u, \chi)$ of
$(P)$, we have $ u\in C^{0,1}_{loc}(\Omega\cup T).$
\end{theorem}

\vs 0,5cm\n The proof of Theorem 1.2 is based on the
next lemma. The rest of the proof will be omitted, since
it can be easily obtained using arguments similar to those in the
proof of Theorem 2.1 \cite{[CL6]} and taking into account the
remark in the proof of Theorem 1.1.

\vs 0,5cm

\begin{lemma}\label{l2.1} Let $S_0$ be an open connected subset of $T$ such that
$S_0\subset\subset T$. Then there exists a positive constant $C$
depending only on $n$, $a$, $M$, $\overline{h}$, $\delta(\Omega)$
  and $R_0$ such that
$$ u(x) \, \leqslant  \,C \, |x-x_0|\qquad\forall x\in \Omega \quad\forall x\in S_0.$$
\end{lemma}

\vspace {0,5cm} \n \emph{Proof.} Let $x_0\in S_0$,  $x_1=x_0+R_0\nu$, where $\nu$ is
the outward unit normal vector to $\partial\Omega$ at $x_0$ such that
$B_{R_0}(x_1)\cap\partial\Omega=\{x_0\}$. Then
we consider the function $v(x)=\vartheta(d(x))$
where $d$ and $\vartheta$ are defined by

$$d(x)=|x-x_1|-R_0,\quad
\vartheta(t)=\int_0^t a^{-1}\Big(\Big(a\Big({M\over
R_0}\Big)+{{\overline{h}R_0}\over{n-1}}\Big)e^{ {{n-1}\over
{R_0}}(D-s)}-{{\overline{h}R_0}\over{n-1}}\Big) ds\quad \text{ and }\quad D=\delta(\Omega).$$

\n Then it is easy to verify the following properties of the function $\vartheta$:
\begin{eqnarray*}
   && \vartheta(0)=0,~ \vartheta(R_0)\geqslant M\\
   && \displaystyle{\vartheta'(t)=a^{-1}\Big(\Big(a\Big({M\over
R_0}\Big)+{{\overline{h}R_0}\over{n-1}}\Big)e^{ {{n-1}\over
{R_0}}(D-t)}-{{\overline{h}R_0}\over{n-1}}\Big)>0\qquad
\forall t\in[0,D]}\\
 && \displaystyle{\vartheta'(D)={M\over R_0}\leqslant \vartheta'(t)\leqslant \vartheta'(0)=a^{-1}
\Big(\Big(a\Big({M\over
R_0}\Big)+{{\overline{h}R_0}\over{n-1}}\Big)e^{ {{n-1}\over
{R_0}}D}-{{\overline{h}R_0}\over{n-1}}\Big)\qquad \forall
t\in[0,D]}
\\
 && \displaystyle{a(\vartheta'(t))\vartheta''(t)+
{{n-1}\over {R_0}}a(\vartheta'(t))+\overline{h}=0 \qquad \forall
t\in[0,D]}.
\end{eqnarray*}

\n We also have

\begin{eqnarray}
&& {{\partial v}\over{\partial x_i}} = \vartheta'(d(x)){{\partial
d}\over{\partial
x_i}}=\vartheta'(d(x)){{x_i-x_{1i}}\over{|x-x_1|}} \nonumber\\
&&{{a(|\nabla v|)}\over{|\nabla v|}}\nabla v=a(\vartheta'(d(x)))\nabla d(x)\nonumber\\
&& {{\partial^2 d}\over{\partial x_i^2}} =
{1\over{|x-x_1|}}-{{(x_i-x_{1i})^2}\over{|x-x_1|^3}}\nonumber\\
&&{{\partial }\over{\partial x_i}}\Big({{a(|\nabla
v|)}\over{|\nabla v|}}{{\partial v}\over{\partial
x_i}}\Big)=a^\prime(\vartheta'(d(x)))\vartheta''(d(x))\Big({{\partial
d}\over{\partial x_i}}\Big)^2+a(\vartheta'(d(x))){{\partial^2
d}\over{\partial x_i^2}}\nonumber\\
&&  \Delta_A v=a^\prime(\vartheta'(d(x)))\vartheta''(d(x))+ {{n-1}\over
{|x-x_1|}}a(\vartheta'(d(x))).\nonumber
\end{eqnarray}

\n Therefore, since $|x-x_1| > R_0$ for all $x$ in $\Omega$, we obtain

\begin{equation}\label{1.8}
\Delta_A v + div(H) \leqslant 0\quad\text{in}\quad \Omega.
\end{equation}

\vs 0,2cm \n Moreover,  we have
\begin{equation}\label{1.9}
u(x)\leqslant v(x)\qquad \text{for all} \quad x\in \partial\Omega.
\end{equation}
\n Indeed, first we have for $x\in T$, $u(x)=0\leqslant
v(x)$. Next, for $x\in\partial\Omega\setminus T$, we have
 $|x-x_0|\leqslant |x-x_1|+|x_1-x_0|=|x-x_1|+R_0$
which  leads to $|x-x_1|\geqslant |x-x_0|-R_0\geqslant
d(S_0,\partial\Omega\setminus T)-R_0=d_0-R_0>3R_0-R_0=2R_0$. This
leads to $ v(x)\geqslant\vartheta(R_0) \geqslant M\geqslant u$ on
$\partial\Omega\setminus T$.

\vs 0,2cm\n Now thanks to (1.9), we have $(u-v)^+\in
W_0^{1,A}(\Omega)$. Using this function in $(P)ii)$ and in (1.8), we
obtain

\begin{eqnarray}\label{1.10-11}
&& \displaystyle{\int_{\Omega}}{{a(|\nabla u|)}\over{|\nabla
u|}}\nabla u  .\nabla(v-u)^+ \,=\,
-\int_{\Omega} \chi H(x) .\nabla(u-v)^+ dx \\
&&\displaystyle{-\int_{\Omega}} {{a(|\nabla v|)}\over{|\nabla
v|}}\nabla v. \nabla(u-v)^+ dx\leqslant \int_{\Omega} H(x)
.\nabla(u-v)^+ dx.
\end{eqnarray}

\n Taking into account that $\chi=1$ a.e. in $\{u>0\}$ and
adding (1.10) and (1.11), we obtain

\begin{equation*}
\int_{\Omega} \Big( {{a(|\nabla u|)}\over{|\nabla u|}}\nabla
u-{{a(|\nabla v|)}\over{|\nabla v|}}\nabla v \Big).\nabla(u-v)^+
dx \leqslant 0
\end{equation*}

\n which leads by (0.2) to $(u-v)^+$ is constant in $\Omega$.
Since $u\leqslant v$ on $\partial\Omega$, we get $u\leqslant v$ in
$\Omega$.

\vs0.3cm\n We conclude that for all $x\in \Omega$ and
$x_0\in S_0$, we have
\begin{eqnarray*}&& u(x)\leqslant v(x)=|v(x)-v(x_0)|\leqslant\sup_{x\in \Omega}|\nabla
v(x)||x-x_0|\\
\qquad&&\leqslant\Big( \sup_{t\in[0,D]}\vartheta'(t)\Big)
|x-x_0|=\vartheta'(0)|x-x_0|=C(n,a,M,\overline{h},D,R_0)|x-x_0|.
\end{eqnarray*}
\qed

\vs 0.5cm
\section{The free boundary}\label{1}

\vs 0.5cm
\n In this section, we assume that the vector function $H$
satisfies  the following assumptions for some positive constants
$\underline{h}$ and $\overline{h}$ :

\begin{eqnarray}\label{2.1-3}
 & |H_1|, |H_2|,..., |H_{n-1}|\leqslant \overline{h},\qquad  0< \underline{h} \leqslant H_n\leqslant  \overline{h}
 \quad\text{ a.e. in }\Omega \\
&  H\in C^{0,1}(\overline{\Omega})\\
&  div(H)\geqslant 0 \quad\text{ a.e. in } \Omega.
\end{eqnarray}

\n By using $\displaystyle{\min\Big({u\over\epsilon},1\Big)\zeta}$
with $\zeta\in {\cal D}(\Omega)$, $\zeta\geqslant 0$ as a test function
for $(P)ii)$ and arguing as in \cite{[CL5]}, one can establish the
following important inequality :

\begin{equation}\label{2.4}
div (\chi H) - \chi(\{u>0\}) div ( H) \leqslant 0 \qquad\text{
in }\quad {\cal D}^\prime(\Omega).
\end{equation}

\n as a consequence of (2.4), we will derive a weak monotonicity of the function  $\chi$,
that will help us to express the free boundary locally as graphs of a family of
functions. More precisely, we consider the following differential equation

\[\displaystyle{(E(\omega,h))\,\,\left\{\begin{array}
{r@{\quad=\quad}l}
X^\prime (t,\omega,h)& H(X(t,\omega,h))\\
 X(0,\omega,h) & (\omega,h)\end{array}\right.}\]

\n where $h\in\pi_{x_n}(\Omega)$ and $\omega\in\pi_{x'}(\Omega\cap \{x_n=h\})$,
$x'=(x_1,...,x_n)$, $\pi_{x'}$ and $\pi_{x_n}$ are respectively the orthogonal projections on
the hyperplane $\{x_n=0\}$ and the $x_n$-axis. We shall denote by $X(.,\omega,h)$ the  maximal solution
of $E(\omega,h)$ defined on the interval $(\alpha_{-} (\omega,h), \alpha_{+}(\omega,h))$ and continuous
in the open set
$$\{ (t,\omega,h)/ \quad
 \alpha_{-} (\omega,h)<t< \alpha_{+}(\omega,h), \,
 h\in\pi_{x_n}(\Omega), \,
 \omega\in \pi_{x'}(\Omega\cap \{x_n=h\})\}.$$

\n We deduce from (2.1) that we have
$$|X(t_1,\omega,h)-X(t_2,\omega,h)|\leqslant \overline{h}|t_1-t_2|~~\forall t_1, t_2\in (\alpha_{-} (\omega,h), \alpha_{+}(\omega,h)).$$
It follows that the limits $\displaystyle{\lim_{t\rightarrow \alpha_{-} (\omega,h)^+}X(t,\omega,h)}$ and
$\displaystyle{\lim_{t\rightarrow \alpha_{+} (\omega,h)^-}X(t,\omega,h)}$ both exit,
which we shall denote respectively by
$X(\alpha_{-} (\omega,h),\omega,h)$ and $X(\alpha_{+} (\omega,h),\omega,h)$, and
observe that we have necessarily $X(\alpha_{-} (\omega,h),\omega,h)\in \partial\Omega\cap
\{x_n<h\}$ and $X( \alpha_{+}(\omega,h),\omega,h)\in\partial\Omega\cap \{x_n>h\}$.

\n For simplicity, we will drop the dependence on $h$ in the sequel.
We shall also denote the orbit of $X(.,\omega)$
by $\gamma(\omega)$.

\vs0.5cm
\n Now, we recall for the reader the following technical properties
and definitions depending only on $H$ which  were established in \cite{[CL5]} :

$\bullet$  $\alpha_{+}$ and $\alpha_{-}$ are uniformly bounded.

\vs0.2cm $\bullet$ For each $h\in\pi_{x_n}(\Omega)$, we define
the set
$$D_h=\{ (t,\omega)\,/\, \omega\in \pi_{x'}(\Omega\cap\{x_n=h\}), \, t\in
 (\alpha_{-}(\omega), \alpha_{+}(\omega))\}$$
 and consider the mapping
 \begin{eqnarray*}
  &&T_h\, :  \,{D}_h \longrightarrow  T_h({D}_h)\\
   &&\qquad(t,\omega)\longmapsto  T_h(t,\omega)=(T_h^1,...,T_h^n)(t,\omega)=X(t,\omega).
 \end{eqnarray*}

\n Clearly each $(x',x_n)\in\Omega$ can be written as
 $(x',x_n)=X(0,\omega)=T_h(0,\omega)$ with $\omega=x'$ and $h=x_n$. So
 $\displaystyle{ \Omega=\bigsqcup_{h\in\pi_{x_n}(\Omega)}T_h(D_h)}$.

\vs0.2cm $\bullet$  $ T_h$ is one to one, $ T_h$  and
$T_h^{-1} $ are $C^{0,1}$.

\vs 0,2cm $\bullet$  The Jacobian  matrix  of the
mapping $T_h$, denoted by $\mathcal{J}T_h$, is in $ L^{\infty}(D_h)$ and determinant of
$\mathcal{J}T_h$,  denoted by $Y_h(t,\omega)$,  satisfies :

\vs 0.2cm\n i) \quad $\displaystyle{\partial Y_h\over
\partial t} (t,\omega) = Y_h(t,\omega) (div H)(X(t,\omega))
\qquad\mbox{a.e. in } D_h.$

\vs 0.2cm \n ii) \quad $\displaystyle{ Y_h(t,\omega)=
-H_n(\omega,h) \exp\big( \int_0^t  (div H)(X(s,\omega)) ds\big)}
\qquad\mbox{a.e. in } D_h.$

\vs 0.2cm \n iii) \quad  $\underline{h} \leqslant
-Y_h(t, \omega) \leqslant C {\overline{h}},\quad C>0,
\qquad\mbox{a.e. in } D_h.$

\vs 0.5cm \n Using (2.4) and arguing as in the proof of Theorem 2.1 of
\cite{[CL5]}, we can establish the following monotonicity of
$\chi$

\begin{equation}\label{1}
  {\partial \over\partial t} \big(\chi o T_h \big)\leqslant 0\qquad\text{ in }
{\cal D}^\prime (D_h).
\end{equation}

\vs 0.5cm\n Property (2.5) means that $\chi$ decreases along the orbits
$\gamma(\omega)$ of the differential equation $(E(w,h))$. The consequence on $u$
is materialized in the next key
theorem which is the main idea in the parametrization of the free boundary.

\begin{theorem}\label{t2.1.} Let $(u,\chi)$ be a solution of
$(P)$ and $X_0=(x_0', x_{n0})= T_h(t_0,\omega_0)\in T_h(D_h)$.

\vs 0,2cm \n
 i)\quad If $\quad u(X_0)=uoT_h(t_0,\omega_0) > 0,\quad$ then there exists $\epsilon >0$
  such that
$$uoT_h(t,\omega) > 0 \qquad \forall
 (t,\omega)\in C_\epsilon=\{ (t,\omega)\in D_h\,/ \,\, \vert \omega-
\omega_0\vert <\epsilon , \,\, t <t_0 +\epsilon \}.$$
\vs 0,2cm
\n ii)\quad If $\quad u(X_0)=uoT_h(t_0,\omega_0)=0,\quad $
then $\quad uoT_h(t, \omega_0) =0 \qquad \forall t\geqslant
t_0.$\end{theorem}

\vs0.5cm \n To prove Theorem 2.1, we need the following strong maximum principle:

\begin{lemma}\label{l2.1}

\n Let $u\in W^{1,A}(U)\cap C^1(U)\cap C^0(\overline{U})$ such that $ u
\geqslant 0$ in $U$ and $\Delta_A u\leqslant 0$ in $U$.
Then $u\equiv 0$ in $U$ or $u>0$ in $U$.
\end{lemma}

\vs 0,3cm \n The proof of Lemma 2.1 follows from the next lemma as in
\cite{[E]} p. 333.

\begin{lemma}\label{l2.2}

\n Let $u\in W^{1,A}(U)\cap C^1(\overline{U})$ such that
$\Delta_A u \leqslant 0$ in $U$ and $u(x_0)<u(x)$ for all
$x\in U$, where
$x_0\in \partial B_R(x_1)$ and $B_R(x_1)\subset U$.

\n Then the outer normal derivative of $u$ at $x_0$,
satisfies ${{\partial u}\over{\partial\nu}}(x_0)<0.$
\end{lemma}

\vs 0,5cm \n \emph{Proof.} We consider the standard function $v$ defined  by
 $$v(x) =
 e^{-\alpha r^{2}}-e^{-\alpha R^{2}} \qquad \hbox{ in } D= B_R(x_1)\setminus \overline{B_{R/2} (x_1)}$$
where
$$r= |x-x_1|\in (R/2, R),\qquad  \alpha={\kappa\over{(R/2)^2}}={{4\kappa}\over{R^2} },$$

\n and $\kappa$ is a positive parameter such that $\displaystyle{ {1\over2}<\kappa<2\left(1+{{n-2}\over{a_0}}\right) }$.

\n Let $\epsilon= \displaystyle{ { \displaystyle{\min_{\partial B_{R/2}(x_1)}  }(u-u(x_0))}\over{ \displaystyle{\max_{\partial B_{R/2}(x_1)}  } v}}
>0$, and observe that

\begin{equation} \label{e2.6}
\epsilon v \leq u-u(x_0)  \qquad \hbox{ on } \partial
D.\end{equation}

\vs 0,2cm\n To establish the Lemma, we will compare $u-u(x_0)$ with respect to
$\epsilon v$. We claim that

\begin{equation}\label{e2.7}
  \Delta_A (\epsilon v)\geqslant {a(\epsilon|\nabla  v|)\over r}\geqslant0\qquad \hbox{ in }  D.
\end{equation}

\n Indeed, we first observe from (1.2) that

\begin{equation}\label{e2.8}
 \displaystyle{ \Delta_A (\epsilon v)= {{a(\epsilon|\nabla v|)}\over{|\nabla v|^3}}\Big\{|\nabla v|^2
  \Delta v+\Big({{a'(\epsilon|\nabla v|)}\over{a(\epsilon|\nabla v|)}}|\epsilon\nabla
  v|-1\Big)\sum_{i,j}{{\partial v}\over{\partial x_i}}{{\partial v}\over{\partial x_j}}
   {{\partial^2 v}\over{\partial
x_i\partial x_j}}\Big\}}.
\end{equation}
\n Moreover we have in $D$

\begin{eqnarray*}
&& \nabla v=-2\alpha k e^{-\alpha r^{2}}(x-x_1), \qquad |\nabla
v|=2\alpha r e^{-\alpha r^{2}},
 \qquad \Delta v=-2\alpha e^{-\alpha r^{2}}(n-2\alpha r^2),\\
&& {{\partial^2 v}\over{\partial x_i\partial x_j}}=-2\alpha
e^{-\alpha r^{2}}\Big(\delta_{ij}-2\alpha(x_i-x_{1i})(x_j-x_{1j})\Big),\\
&& \sum_{i,j}{{\partial v}\over{\partial x_i}}{{\partial
v}\over{\partial x_j}} {{\partial^2 v}\over{\partial x_i\partial
x_j}} =-(2\alpha)^3 r^2e^{-3\alpha r^2}(1-2\alpha r^2).
\end{eqnarray*}

\n Taking into account the fact that $\displaystyle{
1-2\alpha r^2=1-2{\kappa\over r^2}r^2\leqslant
1-2{\kappa\over (R/2)^2}\Big({R\over 2}\Big)^2=1-2\kappa<0}$,
we get by substituting the above formulas into (2.8)

\begin{eqnarray*}
&& \displaystyle{ \Delta_A (\epsilon v)= (2\alpha)^3
r^2e^{-3\alpha r^2}{{a(\epsilon|\nabla v|)}\over{|\nabla
 v|^3}}\Big\{n-1
 +{{a'(\epsilon|\nabla v|)}\over{a(\epsilon|\nabla v|)}}|\epsilon\nabla
  v|(1-2\alpha r^2)\Big\}}\nonumber\\
 && \qquad \geqslant (2\alpha)^3
r^2e^{-3\alpha r^2}{{a(\epsilon|\nabla v|)}\over{|\nabla
 v|^3}}\Big\{n-1
 +a_0(1-{\kappa\over 2})\Big\}\qquad\hbox{ by } (0.1)\nonumber\\
 &&\qquad  ={{a(|\nabla v|)}\over{r}}\Big(n-1
 +a_0(1-{\kappa\over 2})\Big)\geqslant{{a(\epsilon|\nabla v|)}\over{r}}.
\end{eqnarray*}

\n Hence (2.7) holds.

\vs 0.3cm\n Using (2.7) and the fact that $\Delta_A u\leqslant 0$,
we obtain
\begin{equation}\label{e2.9}
\Delta_A (\epsilon v)\geqslant 0\geqslant\Delta_A (u-u(x_0))  \qquad \hbox{ in }D.
\end{equation}
Now taking into account (2.6) and (2.9), and the weak maximum principle for
the $A$-Laplacian, we get

\begin{equation}\label{e2.10}
\epsilon v \leqslant u-u(x_0)  \qquad \hbox{ in }D.
\end{equation}

\n To conclude, let $\nu$ be the exterior unit normal at $x_0$. We
infer from (2.10) for $t$ positive and small enough so that $x_0-t\nu\in D$
$$ { u(x_0-t\nu)- u(x_0)\over t} \geqslant \epsilon { v(x_0-t\nu)- v(x_0)\over
t}.$$

\n Letting $t\longrightarrow 0$, we obtain
\begin{eqnarray*}
&&-{\partial u \over \partial\nu} (x_0)  \geqslant \epsilon .\Big(-{\partial v \over \partial\nu}
(x_0)\Big)\\
&&{\partial u \over \partial\nu} (x_0)  \leqslant \epsilon {\partial v \over \partial\nu}(x_0)= \epsilon.(-2\alpha R
e^{-\alpha R^2}) <0.
\end{eqnarray*}
\qed

\vs 0.5cm \n \emph{Proof of Theorem 2.1.} It is enough to verify $i)$. By
continuity, there exists $\epsilon
>0$ such that
$$uoT_h(t,\omega) > 0 \qquad \forall
 (t,\omega)\in (t_0 -\epsilon, t_0+\epsilon)\times B_{\epsilon}(\omega_0)
= Q_\epsilon.$$

\n By $(P)i)$, we have $\chi oT_h(t,\omega) = 1 $
for a.e. $(t,\omega)\in Q_\epsilon$. Using $(2.5)$ and since $\chi
oT_h\leqslant 1$,  we get $\chi oT_h = 1$ a.e. in $C_\epsilon$,
 i.e. $\chi = 1$ a.e. in $T_h(C_\epsilon)$.

\n From $(P)ii)$ and (2.3),  we get $\triangle_A u=- div(H)\leqslant 0$ in
 ${\cal D}^\prime ({ T}_h(C_\epsilon))$.  Given that $u\geqslant 0$ in $\Omega$ and $u>0$ in
 ${ T}_h(Q_\epsilon)\subset { T}_h(C_\epsilon)$, we conclude by Lemma 2.1,
that $u>0$ in ${ T}_h(C_\epsilon)$.\qed

\vs 0,5cm \n \begin{remark}Thanks to Theorem 2.1, we can then define for each
$h\in \pi_{x_n}(\Omega)$, the following function $\phi_h$ on
$\pi_{x'}(\Omega\cap \{x_n=h\})$ by

\begin{equation}\phi_h(\omega) =\begin{cases}
&\sup\left\{t\,~:~ \, (t,\omega)\in D_h, \quad u o T_h(t,\omega) > 0\right\}
 \\
&\hskip 2cm \text{ if this set is not empty }\\
 & \alpha_-(\omega)\hskip 2cm
\text{ otherwise.}\end{cases}\end{equation}

\n One can easily check as in \cite{[CL1]} that $\phi_h$ is lower
semi-continuous at each $\omega \in\pi_{x'}(\Omega\cap \{x_n=h\})$  such
that $T_h(\phi_h(\omega), \omega) \in\Omega$ and that
$\left\{uoT_h(t,\omega)>0\right\}\cap D_h=\left\{t<\phi_h(\omega)\right\}$.

\end{remark}

\vspace{0,5cm} \n\emph{Acknowledgments } The authors are
grateful for the facilities and excellent research conditions at
the Fields Institute where part of this research was carried out.

\end{document}